\documentclass[a4paper, 12pt]{article} 

\usepackage{amsfonts}
\usepackage{amsmath}
\usepackage{enumerate}
\usepackage[francais]{babel}
\usepackage[applemac]{inputenc}
\usepackage{graphicx}
\usepackage{geometry}
\usepackage{makeidx} 

\newtheorem{lem}{Lemma}
\newtheorem{prop}{Proposition}

\newtheorem{corol}[prop]{Corollary}
\newtheorem{thm}[prop]{Theorem}

\newtheorem{rmq}{Remark}

\newcommand{\R}{\mathbb{R}}
\newcommand{\N}{\mathbb{N}}

\newcommand{\E}{\mathbb{E}}

\newcommand{\C}{\mathbb{C}}

\newcommand{\1}{\mathbf{1}}
\usepackage[bf]{caption}

\makeindex 

\title{Stochastic representation of tau functions of Korteweg--de Vries equation}
\author{Mich\`ele Thieullen\thanks{
Laboratoire de Probabilit\'es et Mod\`eles Al\'eatoires, Universit\'e Pierre et Marie Curie,  4, Place Jussieu, 75252 Paris cedex 05, France; 
michele.thieullen@upmc.fr}\\
\and 
Alexis Vigot\thanks{Laboratoire de Probabilit\'es et Mod\`eles Al\'eatoires, Universit\'e Pierre et Marie Curie,  4, Place Jussieu, 75252 Paris cedex 05, France; 
alexis.vigot@upmc.fr}
}

 \date{24 april 2017}

\begin{document}
\maketitle
\noindent {\bf Abstract :} In this paper we express the tau functions considered by P\"oppe in \cite{poppe2} for the Korteweg de Vries (KdV) equation, as Laplace transforms of iterated Skorohod integrals. Our main tool is the notion of Fredholm determinant of an integral operator. Our result extends the paper \cite{taniguchi} by Ikeda and Taniguchi who obtained a stochastic representation of tau functions for the $N$-soliton solutions of KdV as the Laplace transform of a quadratic functional of $N$ independent Ornstein-Uhlenbeck processes. Our general result goes beyond the $N$-soliton case and enables us to consider a non soliton solution of KdV associated to a Gaussian process with Cauchy covariance function.  



\vspace{1cm}


\section{Introduction}
The solutions of different families of partial differential equations (pdes) can be expressed as the mean of  functionals of a stochastic process by the Feynman-Kac formula (cf. \cite{dynkin},  \cite{pardoux}, \cite{yor}). In this paper, we are interested in KdV equation for which no Feynman-Kac formula exists to represent its solutions. In this paper we prove a stochastic representation for the {\it tau functions} obtained by P\"oppe in \cite{poppe2} for KdV which form a class of fairly general tau functions. Indeed for a pde such as KdV, solutions can be generated using tau functions which are such that derivatives of their logarithm provide solutions of KdV. Tau functions are also available for the KdV (resp. KP) hierarchy that consists in a system of pdes in infinitely many variables which contains KdV (resp. KP) equation (cf. \cite{miwa}). Each equation in these hierarchies admits solutions that can be written as derivatives of the logarithm of a tau function. We will be more precise at the end of this section. We prove that the tau functions of \cite{poppe2} admit a stochastic representation as Laplace transforms of some iterated Skorohod integrals. To illustrate the result we consider a particular situation where the iterated Skorohod integrals can be written as integrals of the square of a Gaussian process with respect to some measure. This measure is discrete in the $N$-soliton case. 

\bigskip

The stochastic representation of tau functions have been studied by Ikeda and Taniguchi in \cite{taniguchi} and Aihara, Akahori, Fujii and Nitta in \cite{tau}.  The case of $N$-soliton solutions of KdV hierarchy is addressed in  \cite{taniguchi} and has been extended  later in \cite{tau} to the $N$-soliton solutions of KP hierarchy. In both papers, the tau function is obtained as the Laplace transform of a quadratic functional of a Gaussian process. Our result extends \cite{taniguchi} to solutions of KdV that are not necessarily $N$-solitons. One aspect of our extension is the use of the Skorohod integral which generalises the Ito integral employed in \cite{taniguchi}.


\bigskip

The key tool to obtain our result is the Fredholm determinant of an integral operator. Fredholm determinants are an extension of the determinant of a matrix to a larger class of operators (cf. \cite{trace}). Indeed in  \cite{poppe2} (see also \cite{poppe}) tau functions for KdV and KP hierarchies can be expressed as the Fredholm determinants of some integral operators as a consequence of the dressing method introduced by Zakharov and Shabat in \cite{zakharov}. Note that in Probability Theory, Fredholm determinants appear in the change of measure formulas on Wiener spaces (cf. for instance to \cite{cameron2} and \cite{ramer}) and also in Laplace transforms (cf. \cite{grasselli}, \cite{privault}). Our first step will be to identify the Fredholm determinant of an integral operator in the stochastic representation of \cite{taniguchi} and to provide the Wiener chaos decomposition of the underlying process. The results of the present paper can be extended to KP and KdV hierarchy.



 of the present paper can be extended to KP and KdV hierarchy.\\

The KdV equation is one of the simplest and most useful nonlinear equations admitting solitary waves, and it represents in particular the time evolution of waves in shallow waters. In the present manuscript we consider the two usual versions of KdV namely
\begin{equation}\label{KdV}
u_t-6uu_x+u_{xxx}=0.
\end{equation}
and
\begin{equation}\label{eq:kdv4}
u_t=\frac{3}{2}uu_x+\frac{1}{4}u_{xxx}.
\end{equation}
Indeed the first one is the most encountered in the litterature and the second one (\ref{eq:kdv4}) is the form that appears when one is interested in KdV hierarchy. We work with (\ref{eq:kdv4}) in Section 2 and (\ref{KdV}) in Section 3. One can go from one form to the other by a linear change of variables.\\

The notion of tau function comes from the Inverse Scattering Transform (IST) introduced by Gardner, Greene, Kruskal and Miura in \cite{ggkm} and \cite{ggkm2}, which is a method for solving some non-linear pdes such as KdV equation and Kadomtsev--Petviashvili (KP) equation
\begin{equation}\label{eq:KP}
\frac{3}{4}u_{yy}=\partial_x\left(u_t-\frac{3}{2}uu_x-\frac{1}{4}u_{xxx}\right).
\end{equation}
This method applied to KdV equation (\ref{eq:kdv4}) provides solutions called $N$-solitons,
\begin{equation}\label{eq:solkdv}
u(x,t):=2\frac{d^2}{dx^2}\log \det (I+G(x,t))
\end{equation}
where $G$ is a $N\times N$-matrix defined by 
\begin{equation}\label{eq:200}
G_{ij}(x,t)=\frac{\sqrt{m_im_j}}{\eta_i+\eta_j}e^{-(\eta_i+\eta_j)x-(\eta_i^3+\eta_j^3)t}
\end{equation}
and $m_i,\eta_i>0$, $1\leq i\leq N$ are constants called scattering data. For more details we refer the reader to \cite{drazin}. In particular identity (\ref{eq:solkdv}) states that $\det (I+G(x,t))$ is a tau function for KdV. A tau function for KP hierarchy providing a \emph{$N$-soliton} solution is
\begin{equation}\label{taukp}
\tau(x_1,x_2,\ldots)=\det(I+G(x_1,x_2,\ldots)),
\end{equation}
where $G$ is the square matrix
\begin{equation}\label{eq:244}
G(x_1,x_2,\ldots):=\left(\frac{\sqrt{m_im_j}}{p_i-q_j}e^{-\frac{1}{2}\sum_{l=1}^\infty \{(p_i^l-q_i^l)+(p_j^l-q_j^l)\}x_l}\right)_{1\leq i,j\leq N},
\end{equation}
and $p_i\in\R$, $q_i\in \R$ and $m_i> 0$, $1\leq i\leq N$ are constants. The tau function for the $N$-soliton of KdV hierarchy is obtained from (\ref{eq:244}) by taking $p_i=\eta_i$ and $q_i=-\eta_i$, $1\leq i\leq N$ for positive constants $\eta_i$.\\

In Section 2 we recall definitions of Fredholm determinant and Carleman-Fredholm determinant. In Section 3, we present the result from \cite{taniguchi} and write it as the Fredholm determinant of an integral operator. In Section 4, result from \cite{tau} are presented. In Section 5, we give our stochastic representation of a tau function for KdV and retrieve the $N$-soliton solutions as a particular case. In Appendix a Fubini theorem is given for Skorohod integrals. 

\section{Reminders on Fredholm determinants}
In this section, we recall some definitions on Fredholm determinants and we refer to \cite{trace} for more details.\\

\noindent Fredholm in \cite{fredholm} showed that for an integral operator 
\begin{equation}\label{eq:300}
\mathcal{A}f(x):=\int_a^b K(x,y)f(y)dy
\end{equation}
 with continuous kernel $K$, the integral equation $(I+z\mathcal{A})\phi=\psi,$ where $\psi$ is a given continuous function, admits a unique continuous solution $\phi$ if and only if $D(z)\neq 0$ with $D$ the entire function
\begin{equation}\label{eq:217}
D(z):=1+\sum_{\mu=1}^\infty\frac{z^\mu}{\mu!}\int_a^b\ldots \int_a^b\begin{vmatrix} K(s_1,s_1) & \cdots & K(s_1,s_\mu) \\ \vdots & \ddots & \vdots\\ K(s_\mu,s_1) & \cdots & K(s_\mu,s_\mu)\end{vmatrix}ds_1\ldots ds_\mu, \quad z\in \C.
\end{equation}
$D(z)$ is called the Fredholm determinant of $I+z\mathcal{A}$ and is written $\det(I+z\mathcal{A})$. \\

Let $H$ be a complex, separable Hilbert space and let $\mathcal{J}_\infty(H)$ be the set of compact operators on $H$. Let $\mathcal{A}\in \mathcal{J}_\infty(H)$, denoting by $(s_n(\mathcal{A}))_{n=1}^N$, $N\leq\infty$ the eigenvalues of the positive semidefinite operator $|\mathcal{A}|:=(\mathcal{A}^\star\mathcal{A})^{1/2}$, we define a family of embedded subsets of $ \mathcal{J}_\infty(H)$
$$ \mathcal{J}_p(H):=\{\mathcal{A}\in  \mathcal{J}_\infty(H), \sum_{n=1}^Ns_n(\mathcal{A})^p<\infty\},\quad 1\leq p<\infty.$$
The operators in $\mathcal{J}_1(H)$ are called \emph{trace class operators} and those in $\mathcal{J}_2(H)$ are called \emph{Hilbert-Schmidt operators}.\\

\noindent Let $\mathcal{A}$ be a trace class operator with eigenvalues $(\lambda_n(\mathcal{A}))_{n=1}^N$, then its trace is defined by
\begin{equation}\label{eq:255}
\text{Tr\ }\mathcal{A}:=\sum_{n=1}^N\lambda_n(\mathcal{A}),
\end{equation}
and the Fredholm determinant $\det(I+z\mathcal{A})$ is defined by the product
\begin{equation}\label{eq:239}
\det(I+z\mathcal{A}):=\prod_{n=1}^N(1+z\lambda_n(\mathcal{A})).
\end{equation}
These two identities make sense because of the assumptions on $\mathcal{A}$. In particular, if the integral operator (\ref{eq:300}) is of trace class its trace is simply
\begin{equation}\label{eq:247}
\text{Tr\ }\mathcal{A}=\int_a^bK(x,x)dx.
\end{equation}
If $\mathcal{A}$ is an Hilbert-Schmidt operator, the right-hand side of (\ref{eq:239}) doesn't necessarily converge (when $N=+\infty$). However, in that case we can define the Carleman-Fredholm determinant 
\begin{equation}\label{eq:256}
{\det}_2(I+z\mathcal{A}):=\prod_{n=1}^N(1+z\lambda_n(\mathcal{A}))e^{-z\lambda_n(\mathcal{A})}.
\end{equation}
From the previous definitions, the Carleman-Fredholm determinant of a trace class operator $\mathcal{A}$ can be expressed w.r.t. its Fredholm determinant and trace as
$${\det}_2(I+z\mathcal{A})=\det(I+z\mathcal{A})e^{-z\text{Tr\ }\mathcal{A}}.$$

\section{Fredholm determinant and the stochastic representation in \cite{taniguchi}}\label{sec:taniguchi2}

The following theorem has been obtained in \cite{taniguchi} for KdV hierarchy and is inspired from a result of Cameron and Martin (cf. \cite{cmsturm}). We present it here for KdV equation.\\

Let $N\in\N^\star$, $a>0$, $p\in \R^N$, $D=diag[p_1,\ldots,p_N]$ where $p_i\neq p_j$ for $i\neq j$ and let $c\in \R^N$ be a positive vector. Let
\begin{equation}\label{eq:243}
E(a):=D^2+a^2c\otimes c.
\end{equation}
Then $E(a)$ can be written as $E(a)=UR^2U^{-1}$, where $R$ is a diagonal matrix and $U$ an orthogonal matrix. Let $\xi_p$ be an $\R^N$-valued Ornstein-Uhlenbeck process defined as the unique solution of the stochastic differential equation
\begin{equation}\label{eq:tan}
d\xi_p(s)=dW(s)+D\xi_p(s)ds,\quad \xi_p(0)=0,
\end{equation}
where $W$ is a $N$-dimensional Wiener process. We define the process $X_{p,c}$ by
\begin{equation}\label{eq:320}
X_{p,c}(s):=\langle c,\xi_p(s)\rangle, \quad s\geq 0.
\end{equation}

\begin{thm}(cf. \cite{taniguchi}) \label{thm:taniguchi}
Let $\beta_{a,x,t}$ be defined by $\beta_{a,x,t}(y):=-((\partial_z\phi_a)\phi_a^{-1})(x-y,t)$, where
\begin{equation}
\label{eq:141}
\phi_a(z,t):=U\{\cosh(\zeta(z,t))-\sinh(\zeta(z,t))R^{-1}U^{-1}DU\}U^{-1}
\end{equation}
and $\zeta(z,t)=zR+tR^3$. 
Define $I_{p,c,a}$ by
\begin{equation}
\label{eq:taniguchi}
I_{p,c,a}(x,t):=\E\left[e^{-\frac{a^2}{2}\int_0^xX_{p,c}(y)^2dy+\frac{1}{2}\langle(\beta_{a,x,t}(x)-D)\xi_p(x),\xi_p(x)\rangle}\right].
\end{equation}
Then
$$u(x,t):=-4\ \partial_x^2\log I_{p,c,a}(x,t)$$
satisfies KdV equation (\ref{eq:kdv4}).
\end{thm}

\noindent The key argument to prove this theorem is the identity
\begin{equation}\label{eq:tan2}
I_{p,c,a}(x,t)=\left({\det \phi_a(0,t)}/{\det \phi_a(x,t)}\right)^{1/2}\cdot e^{-\frac{x}{2}\text{Tr\ }D},
\end{equation}
which is a consequence of Girsanov Theorem and the fact that $\phi_a$ satisfies $\partial_{xx}\phi_a-E(a)\phi_a=0.$\\

\noindent In this section, we stress the fact that $\left({\det \phi_a(0,t)}/{\det \phi_a(x,t)}\right)^{1/2}$ coincides with the Fredholm determinant of some integral operator that we make explicit below. 
Let us consider the Volterra integral operator defined on the space of continuous functions from $[0,x]$ to $\R^d$ by
\begin{equation}
\label{eq:146}
V(f):=\int_0^\cdot\phi_a'(u)\phi_a^{-1}(u)f(u)du.
\end{equation}

\begin{prop}\label{prop:dettanig}
Let $V$ be defined by (\ref{eq:146}). Then $I_{p,c,a}$ defined in (\ref{eq:taniguchi}) is a Fredholm determinant namely
$$
I_{p,c,a}(x,t)=\det(I-V/2).
$$
\end{prop}

\noindent\textbf{Proof of Proposition \ref{prop:dettanig}.}
The kernel of $V$,
$$K(v,u):=\phi_a'(u)\phi_a(u)^{-1}\1_{0\leq u\leq v\leq x},$$
being continuous, the trace of $V$ is well defined and given by $\text{Tr\ }V=\int_0^x K(s,s)ds.$
From \cite{trace}, the Fredholm determinant of $T:=I-\lambda V$, $\lambda\in \R$ can be obtained from the trace of $V$ and of its powers as follows
$$\det(I-\lambda V)=\exp\left(-\sum_{n=1}^\infty \frac{\text{Tr\ }V^n}{n}\lambda^n\right).$$
Moreover $\text{Tr\ }V^n=0$ for all $n>1$ as a consequence of Lemma \ref{lem:karambal} below. Using Jacobi's formula, we obtain
\begin{align*}
\det(I-\lambda V)&=\exp \left(-\lambda \text{Tr\ }\int_0^x K(s,s)ds\right)\\
&=\exp \left(-\lambda\int_0^x\text{Tr\ }(\phi_a'(s)\phi_a^{-1}(s))ds\right)=\left(\frac{\det(\phi_a(0))}{\det(\phi_a(x))}\right)^\lambda.
\end{align*}
We conclude by choosing $\lambda=1/2$ and thanks to identity (\ref{eq:tan2}).
\hfill \framebox[0.6em]{}\\

\begin{lem}(cf. \cite{karambal})\label{lem:karambal}
Let $T_1$, $T_2,\ldots, T_\ell$ be Hilbert-Schmidt operators with respective kernels $G_1$, $G_2$, \ldots, $G_\ell$ in $L^2(\R^2;\C^{N\times N})$. Then we have for $\ell\geq 2$
$$\text{Tr\ }(T_1T_2\cdots T_\ell)=\text{Tr\ }\int_{\R^\ell}G_1(s_1,s_2)G_2(s_2,s_3)\cdots G_\ell(s_\ell,s_1)ds_\ell\cdots ds_1.$$
\end{lem}

\subsection{Fredholm determinant and scattering data}

\label{sec:wienerchaos}
We have seen in the previous section that the stochastic representation of a tau function of KdV given in \cite{taniguchi} actually involves the Fredholm determinant of the operator $V$ that we have defined in (\ref{eq:146}). However, this operator depends on $\phi_a$. It would be preferable to obtain a direct correspondance between the integral operator and the Ornstein-Uhlenbeck processes (cf. (\ref{eq:tan})) involved in the stochastic representation of Theorem \ref{thm:taniguchi}. In the present section we exhibit an integral operator whose kernel depends explicitly on the coefficients of these processes and such that its Fredholm determinant coincides with $I_{p,c,a}$.\\ 

Let $\Delta:=[0,x]\times \{1,\ldots,N\}$ equipped with $\mu=\lambda \otimes \nu$ the tensor product of the Lebesgue measure on $[0,x]$ with the counting measure on $\{1,\ldots,N\}$.

\begin{thm}\label{prop:9}Let $\mathcal{C}$ denote the integral operator $\mathcal{C}\psi(\cdot):=\int_\Delta f_2(\cdot,\eta)\psi(\eta)d\mu(\eta),$ where
\begin{multline*}
f_2((v,j),(u,i)):=a^2\frac{c_ic_j}{p_i+p_j}\left(e^{(p_i+p_j)x}-e^{(p_i+p_j)u\vee v}\right)e^{-p_iu}e^{-p_jv}\\
-(\beta_{a,x,t}(x)-D)_{i,j}e^{(p_i+p_j)x-up_i-vp_j},
\end{multline*}
and $u \vee v$ denotes  $\max(u,v)$. Then
\begin{equation}
I_{p,c,a}(x,t)=\left(\det(I+\mathcal{C})\right)^{-1/2}.
\end{equation}
\end{thm}
\bigskip

In the following we will use the notations of \cite{grasselli} and \cite{nualart} for the stochastic integration $I_1(f_1)$ of a function $f_1\in L^2(\Delta,\mu)$
$$I_1(f_1)=\int_\Delta f_1(\eta)dW_\eta=\sum_{i=1}^N\int_0^xf_1(u,i)dW_u^i,$$
and the double stochastic integral $I_2(f_2)$ of a symmetric function $f_2\in L^2(\Delta^2,\mu^2)$
$$\frac{1}{2}I_2(f_2)=\int_{\Delta_2} f_2(\eta_1,\eta_2)dW_{\eta_1}dW_{\eta_2}=\sum_{i_1,i_2=1}^N\int_0^x\int_0^{s_2}f_2((s_1,i_1),(s_2,i_2))dW_{s_1}^{i_1}dW_{s_2}^{i_2},$$
where $\Delta_2:=\{((s_1,i_1),(s_2,i_2))\in \Delta^2, 0\leq s_1\leq s_2\leq x\}$, $(W^i)_{1\leq i\leq N}$ are $N$ independant Brownian motions.\\

\noindent The proof of Theorem \ref{prop:9} makes use of the following proposition.

\begin{prop}(cf. \cite{grasselli})\label{prop:grasselli} Let $B\in L^2(\Delta,\mu)$ and let $C$ be the kernel of a symmetric Hilbert-Schmidt operator $\mathcal{C}$ on $L^2(\Delta,\mu)$ such that $(I + \mathcal{C})$ has positive spectrum. Let $Y$ be a random variable admitting the following Wiener chaos decomposition
$$Y=\int_\Delta B(\eta)dW_\eta+\int_{\Delta_2} C(\eta_1,\eta_2)dW_{\eta_1}\ dW_{\eta_2}.$$
Then $\E\left[e^{-Y}\right]$ is well defined satisfies
\begin{equation}\label{eq:303}
\E\left[e^{-Y}\right]=[{\det}_2(I+\mathcal{C})]^{-1/2}\exp\left[\frac{1}{2}\int_{\Delta_2} B(\eta_1)(I+\mathcal{C})^{-1}(\eta_1,\eta_2)B(\eta_2)d\eta_1\ d\eta_2\right],
\end{equation}
where ${\det}_2$ is the Carleman-Fredholm determinant and we use the same notation for the operator $(I+\mathcal{C})^{-1}$ and its kernel.\\
\end{prop}

\begin{rmq}\label{rmq:1}
It has been proved in \cite{chiang} and \cite{grasselli} that the existence of the expectation (\ref{eq:303}) is actually equivalent to the assumption that $(I + \mathcal{C})$ has positive spectrum.\\
\end{rmq}

\noindent\textbf{Proof of Theorem \ref{prop:9}.} For simplicity we write in detail the proof for the case $t=0$. The arguments are similar for $t>0$. When $t=0$, $\beta_{a,x,t}(x)=D$ (cf. Section \ref{sec:taniguchi2}). Let $f$ be the $\R^N$-valued function given by $f_i(u):=c_ie^{p_i(y-u)}\1_{[0,y]}(u)$. 
Then, $X_{p,c}$ defined in (\ref{eq:320}) satisfies
$$X_{p,c}(y)=\sum_{i=1}^N c_i\int_0^ye^{p_i(y-u)}dW_u^i=I_1(f).$$
Denoting by $D$ the Malliavin derivative operator, we have for $(u,i)\in \Delta$,
$$D_{(u,i)}X_{p,c}^2(y)=2\ D_{(u,i)}X_{p,c}(y)\cdot X_{p,c}(y)=2\ \1_{[0,y]}(u)c_ie^{p_i(y-u)}X_{p,c}(y).$$
Hence,
\begin{align*}
D_{(u,i)}&\int_0^xX_{p,c}^2(y)dy\\
&=\int_0^x2\ \1_{[0,y]}(u)c_ie^{p_i(y-u)}X_{p,c}(y)dy\\
&=\int_0^x2c_ie^{p_i(y-u)}\1_{[0,y]}(u)\sum_{j=1}^N\int_0^yc_je^{p_j(y-v)}dW_v^jdy\\
&=2\sum_{j=1}^N\int_0^x\frac{c_ic_j}{p_i+p_j}\left(e^{(p_i+p_j)x}-e^{(p_i+p_j)u\vee v}\right)e^{-p_iu}e^{-p_jv}dW_v^j.
\end{align*}
Thus, $D_{(u,i)}\int_0^xX_{p,c}^2(y)dy=2I_1(f_2(\cdot,(u,i)))$ with 
$$f_2((v,j),(u,i)):=\frac{c_ic_j}{p_i+p_j}\left(e^{(p_i+p_j)x}-e^{(p_i+p_j)u\vee v}\right)e^{-p_iu}e^{-p_jv}.$$
Note that $f_2$ is symmetric with respect to its two variables $(u,i)$ and $(v,j)$. We then obtain the following expansion of $\int_0^xX_{p,c}(y)^2dy$ into a finite sum of multiple stochastic integrals of a symmetric function
\begin{equation}
\label{eq:50}
\int_0^xX_{p,c}(y)^2dy=\sum_{i=1}^N\frac{c_i^2}{4p_i^2}\left(e^{2p_ix}-2p_ix-1\right)+I_2(f_2).
\end{equation}
We can apply Proposition \ref{prop:grasselli} to $I_2(f_2)$. Indeed thanks to (\ref{eq:50}) and Theorem \ref{thm:taniguchi} (with $t=0$) we know that $\E\left[e^{-I_2(f_2)}\right]$ is well defined. Remark \ref{rmq:1} ensures that $(I+\mathcal{C})$ has positive spectrum hence (\ref{eq:303}) is valid and we have
$$\E\left[e^{-I_2(f_2)}\right]=\left({\det}_2(I+\mathcal{C})\right)^{-1/2}.$$
This implies
$$\E\left[e^{-\frac{a^2}{2}\int_0^xX_{p,c}(y)^2dy}\right]=\left({\det}_2(I+\mathcal{C})\right)^{-1/2}e^{-\frac{a^2}{2}\sum_i\frac{c_i^2}{4p_i^2}\left(e^{2p_ix}-2p_ix-1\right)}.$$
The determination of the trace
$$\text{Tr}(\mathcal{C})=\sum_i\int_0^xf_2((u,i),(u,i))du=\sum_{i=1}^N\frac{c_i^2}{4p_i^2}\left(e^{2p_ix}-2p_ix-1\right),$$
ends the proof.
\hfill \framebox[0.6em]{}\\



\section{Approach of Aihara et al. in \cite{tau}}\label{sec:aihara}

In \cite{tau}, the tau function of the $n$-soliton solution of KP hierarchy (\ref{eq:244}) is related to the Laplace transform of some generalized stochastic area functional. The method provides also the tau function for the $n$-soliton solution of KdV hierarchy. Hence \cite{tau} is more general that \cite{taniguchi} but also more recent. We present in this section results of \cite{tau}.\\


Let $W^l=(W^{l,1},W^{l,2})$, $l=1,\ldots,n$ be mutually independent two-dimensional Brownian motions starting at the origin. The stochastic area of $W^l$ is
$$S^l:=\int_0^1(W_s^{l,2}dW_s^{l,1}-W_s^{l,1}dW_s^{l,2}).$$
Let $\Lambda:=\text{diag\ }\{\lambda_1,\ldots,\lambda_n\}$, where $\lambda_l$, $l=1,2,\ldots,n$ are positive numbers. Let $C$ be a $n\times n$ real matrix, and $C^\pm$ be its symmetric and skew-symmetric parts $C^\pm=(C\pm C^\star)/2$. Let us set ${\bf W}_t^i=(W_t^{1,i},\ldots,W_t^{n,i})$ for $i=1,2$, and define for any complex number $\sigma$
$$\hat{S}(\sigma):=\sigma\sum_{l=1}^n\lambda_lS^l+\sigma\langle\Lambda^{1/2}C^-\Lambda^{1/2}{\bf W}_1^1,{\bf W}_1^2\rangle-\frac{\sigma^2}{2}\sum_{i=1,2}\langle\Lambda^{1/2}C^+\Lambda^{1/2}{\bf W}_1^i,{\bf W}_1^i\rangle.$$

\begin{thm}(cf. \cite{tau})\label{thm:tau}
If both $\max_l|\lambda_l|$ and $||C^+||$ are sufficiently small, then we have
\begin{equation}
\E [e^{\hat{S}(i)}]=\det(\cosh \Lambda+C\sinh \Lambda)^{-1}.\label{eq:201}
\end{equation}
\end{thm}
This theorem can be proved by using conditional expectation of $e^{\hat{S}(i)}$ conditioned on ${\bf W}_1=(W_1^1,W_1^2)$ and the well known formula (cf. \cite{levy})
$$\E\left[e^{i\xi S_t}|B_t^1=x,B_t^2=y\right]=\frac{\xi t}{\sinh \xi t}e^{(x^2+y^2)(1-\xi t \text{\ coth\ }\xi t)/2t},$$
where $S_t:=\int_0^tB_s^2dB_s^1-\int_0^tB_s^1dB_s^2$ is the stochastic area of a two-dimensional Brownian motion $B$. We refer the reader to \cite{tau} for details.\\

The tau function (\ref{eq:244}) of the $n$-soliton of KP hierarchy is $\tau(x_1,x_2,\ldots)=\det(I+G(x_1,x_2,\ldots))$, where $G$ is given by (\ref{eq:244})
$$G=\left(\frac{\sqrt{m_im_j}}{p_i-q_j}e^{-\frac{1}{2}(\xi_i+\xi_j)}\right)_{1\leq i,j\leq n},$$
and $\xi\in\R^n$ is defined by $\xi_i=\sum_{l=1}^\infty \{(p_i^l-q_i^l)+(p_j^l-q_j^l)\}x_l$.\\
The main of the following theorem is to express $\tau$ with respect to $\E\ [e^{\hat{S}(i)}]$.

\begin{thm}(cf. \cite{tau})\label{thm:1}
Let $P=(1/(p_i-q_j))_{1\leq i,j\leq n}$,and assume that $\min_{i,j}|p_i-q_j|$ is sufficiently large so that $I+P$ is invertible. Let $C:=(I-P)(I+P)^{-1}$ and $\Lambda:=\text{\ diag}\left\{\frac{1}{2}(\xi_1-\log m_1),\ldots,\frac{1}{2}(\xi_n-\log m_n)\right\}$. Then $e^{-\sum\xi_i/2}\E\ [e^{\hat{S}(i)}]^{-1}$ is the tau function for the $n$ soliton solution of KP hierarchy.
\end{thm}

\noindent\textbf{Proof of Theorem \ref{thm:1}.}
Since $G$ given by (\ref{eq:244}) satisfies $G=e^{-\Lambda}Pe^{-\Lambda}$, we have $\tau=\det(I+e^{-\Lambda}Pe^{-\Lambda})=\det(I+Pe^{-2\Lambda})$.
\begin{align*}
\det(\cosh \Lambda+C\sinh \Lambda)&=\det\left(\frac{e^{\Lambda}+e^{-\Lambda}}{2}+C\frac{e^{\Lambda}-e^{-\Lambda}}{2}\right)\\
&=2^{-n}\det((I+C)e^{\Lambda}+(I-C)e^{-\Lambda})\\
&=2^{-n}\det((I+C)e^{\Lambda})\det(I+(I-C)(I+C)^{-1}e^{-2\Lambda})\\
&=2^{-n}\det(I+C)e^{\sum(\xi_i+\log m_i)/2}\det(I+Pe^{-2\Lambda})
\end{align*}
The latter equality holds true since $C=(I-P)(I+P)^{-1}$ is equivalent to $P=(I-C)(I+C)^{-1}$. Using Theorem \ref{thm:tau} ends the proof.
{\hfill \framebox[0.6em]{}}\\

Let $\eta_i$, $i=1,\ldots,n$ be positive constants. Setting $x_1=x$, $x_3=t$, $x_k=0$ for $k\geq 4$ and for $i=1,\ldots, n$, $p_i=\eta_i$ and $q_i=-\eta_i$, (\ref{eq:244}) becomes the matrix (\ref{eq:200}) and $u(x,t):=2\frac{d^2}{dx^2}\log \tau(x,t)$ is solution of KdV equation (\ref{eq:kdv4}). As a corollary of the previous theorem, we obtain

\begin{corol}\label{corol:5}
Let $P=(1/(\eta_i+\eta_j))_{1\leq i,j\leq n}$,and assume that $\min_{i,j}(\eta_i+\eta_j)$ is sufficiently large so that $I+P$ is invertible. Let $C:=(I-P)(I+P)^{-1}$ and $\Lambda:=\text{\ diag}\left\{\eta_1x+\eta_1^3t-\frac{1}{2}\log m_1,\ldots,\eta_nx+\eta_n^3t-\frac{1}{2}\log m_n\right\}$. Then the $n$-soliton solution of KdV equation is $u=-2\partial_x^2\log \E\ [e^{\hat{S}(i)}]$.
\end{corol}

Note the difference between Theorem \ref{thm:1} and Corollary \ref{corol:5}. The theorem provides a tau function valid for a whole family of pdes while the corollary expresses the solution of a given equation in this family.

\section{Hirota's tau function for KdV and its stochastic representation as a Fredholm determinant}

We consider in this section the KdV equation (\ref{KdV}). We have seen that \cite{tau} and \cite{taniguchi} deal with tau functions for soliton solutions of KP and KdV hierarchies. In that case, the tau function is the determinant of a matrix. In \cite{poppe2}  P\"oppe proposed tau functions that are Fredholm determinants and obtained a more general class of solutions. His results rely on the method developed by Hirota who introduced a new formulation for nonlinear PDEs based on a bilinear operator (\ref{eq:321}) (cf. \cite{hirota1} and \cite{hirota2}). In our probabilistic approach of the solutions of KdV we use the Fredholm determinants of \cite{poppe2}. The solutions will be expressed as a Laplace transform of some iterated Skorohod integrals and the solitonic case will be recovered as a particular case.\\

\noindent It is easy to check that if $u$ is a solution of (\ref{KdV}) and $u=-2\partial_x^2\log \tau$, then $\tau$ is solution of
\begin{equation}\label{eq:218}
\tau\tau_{xt}-\tau_x\tau_t+\tau\tau_{xxxx}-4\tau_x\tau_{xxx}+3\tau_{xx}^2=0.
\end{equation}
Hirota's method is to construct solutions to KdV equation by solving (\ref{eq:218}). In order to do so, he introduced (cf. \cite{hirota2}) a bilinear operator $(a,b)\mapsto {\bf{D}}_t^m{\bf{D}}_x^n(a\cdot b)$ defined for two sufficiently smooth functions $a$ and $b$ of two variables  $(x,t)$ and $m, n\in \N$ by 
\begin{equation}\label{eq:321}
{\bf{D}}_t^m{\bf{D}}_x^n(a\cdot b)=\left. \left(\frac{\partial}{\partial t}-\frac{\partial}{\partial t'}\right)^m\left(\frac{\partial}{\partial x}-\frac{\partial}{\partial x'}\right)^na(x,t)b(x',t')\right|_{x'=x,\ t'=t}.
\end{equation}
For instance, when $n=m=1$, 
$${\bf{D}}_t{\bf{D}}_x(a\cdot b)=a_{xt}b+ab_{xt}-a_tb_x-a_xb_t.$$
With this operator the equation (\ref{eq:218}) takes the simple form
\begin{equation}\label{eq:219}
{\bf{D}}_x({\bf{D}}_t+{\bf{D}}_x^3)(\tau\cdot \tau)=0.
\end{equation}

The $N$-soliton solution of KdV is obtained by looking for solutions of (\ref{eq:219}) as power series in a small parameter $\epsilon$
\begin{equation}\label{eq:238}
\tau=1+\epsilon \tau_1+\epsilon^2\tau_2+\cdots.\end{equation}
If $\tau$ satisfies (\ref{eq:219})  and (\ref{eq:238}) and if moreover we choose $\tau_1=\sum_{i=1}^Nc_ie^{\eta_ix-\eta_i^3t}$ for parameters $c_i$, $\eta_i$ then we find that $\tau_i=0$ for $i>N$ and we can solve (\ref{eq:219}) for $\tau_2, \ldots, \tau_N$ iteratively. The solution $\tau$ obtained coincides with the determinant of the $N\times N$ matrix (cf. \cite{drazin})
\begin{equation}\label{eq:240}
A_{nm}(x,t)=\delta_{nm}+c_n^2\frac{e^{-(\kappa_m+\kappa_n)x}}{\kappa_m+\kappa_n}e^{8\kappa_n^3t}.
\end{equation}
For $N=1$, $\tau(x,t)=1+e^{8t-2x}$ solves (\ref{eq:219}) and the corresponding solution of KdV is the $1$-soliton solution $u=-2\partial_x^2\log \tau(x,t)=-2\text{sech}^2(x-4t)$.\\

In the present section we will express tau functions of KdV obtained by Hirota's method as Fredholm determinants of  integral operators. 
 The interval of integration for these operators being unbounded, we define a set of functions and make assumptions on their kernels such that their trace and their determinant are well defined on this set (cf. \cite{glm}).  Let $1/2<\nu\leq 1$ and define
$$C_\nu:=\{\phi\in C^0([0,+\infty),\C), ||\phi||_\nu:=\sup_{s\geq 0}|\phi(s)|(1+s)^\nu<\infty\},$$
where $C^0([0,+\infty),\C)$ denotes the set of continuous functions from $[0,+\infty)$ to $\C$. We define a normed space $LC_\nu$ of Fredholm integral operators 
\begin{equation}
\label{eq:223}
\mathcal{A}\phi(s):=\int_0^{+\infty}A(s,t)\phi(t)dt, \quad \phi\in C_\nu
\end{equation}
 with continuous kernel $A$ satisfying
$$||\mathcal{A}||_\nu:=\sup_{s,t\geq 0}(1+s)^\nu(1+t)^\nu|A(s,t)|<\infty.$$

For all $\mathcal{A}\in LC_\nu$, the Fredholm determinant of $I+\lambda \mathcal{A}$ and the trace of $\mathcal{A}$ are well defined by (\ref{eq:217}) and (\ref{eq:247}) respectively with $b=+\infty$. Let $p(\lambda):=\det(I+\lambda \mathcal{A})$, then $p$ is analytic on $\C$. If $A$ depends on a parameter $x$ and is differentiable with respect to $x$, then $p$ satisfies
\begin{equation}
\label{eq:229}
\partial_x\ p(\lambda)=p\cdot \text{Tr\ }(\partial_x(\lambda \mathcal{A})(I+\lambda \mathcal{A})^{-1}),
\end{equation}
for all $\lambda\in \C$ such that $(1+\lambda \mathcal{A})$ is invertible.

\begin{thm}(cf. \cite{poppe2})\label{thm:poppe}
Let $F$ be solution of the linearized KdV equation
\begin{equation}\label{eq:220}
F_t+8F_{xxx}=0,
\end{equation}
such that $F$ and its derivatives up to order $4$ in $x$ and order $2$ in $t$ are in $C_{2\nu}$. Let $\mathcal{F}_{(x,t)}$ be the Fredholm integral operator defined by
\begin{equation}\label{eq:222}
\mathcal{F}_{(x,t)}\phi(s):=\int_0^{+\infty} F(s+u+2x,t)\phi(u)du.
\end{equation}
Then,
\begin{equation}\label{eq:221}
\tau(x,t;\lambda):=\det(I+\lambda \mathcal{F}_{(x,t)})
\end{equation}
is solution of (\ref{eq:219}) for every $\lambda\in \C$. Moreover, if $\tau$ is nowhere vanishing,
$$u=-2\frac{\partial}{\partial x} \left(\frac{\tau_x}{\tau}\right)$$
solves KdV equation (\ref{KdV}).
\end{thm}

The proof of the previous theorem is interesting since it is rather simple and can be adapted to other nonlinear equations such as KP equation. It makes use of a continuous functional $[\cdot]$ defined on operators $\mathcal{A}$ of the form (\ref{eq:223}) by
$$[\mathcal{A}]:=-A(0,0).$$ %
Let us denote for simplicity $\mathcal{F}_{(x,t)}$ by $\mathcal{F}$. Note that $\mathcal{F}$ belongs to $LC_\nu$. The functional $[\cdot]$ applied to this operator satisfies
$$[\mathcal{F}]=\text{Tr\ }(\partial_x\ \mathcal{F})\quad \text{and}\quad [(I+\mathcal{F})^{-1}\mathcal{F}]=\text{Tr\ }(\partial_x\ \mathcal{F})(I+\mathcal{F})^{-1}.$$
From (\ref{eq:229}) the latter identity corresponds to the derivative with respect to $x$ of $\log \det (I+\mathcal{F})$.\\

In the following, we will express the tau function given by Theorem \ref{thm:poppe} as the Laplace transform of some second order Wiener functionals. Fredholm determinants are related to Laplace transforms of some Wiener process functionals. Indeed, they arise naturally in the transformation of Wiener measure as in papers of Cameron and Martin \cite{cameron1,cameron2}, in Ramer formula \cite{ramer} and Kusuoka's version of Girsanov's theorem.\\
We denote by $\delta$ the Skorohod integral defined as the adjoint of the Malliavin derivative $D$. For properties of this integral we refer the reader to \cite{nualart} . In the proposition below, the notation $\delta(\delta(\phi))$ denotes the iterated Skorohod integral of $\phi$ which is the adjoint of $D^2$ (cf. \cite{nulartzakai}).

\begin{prop}(cf. \cite{privault})
\label{prop:privault}
Let $\phi\in L^2(\R_+^2)$, such that $||\phi||_{L^2(\R_+^2)}<1$. Then
\begin{equation}
\label{eq:152}
\E[e^{-\frac{1}{2}\delta(\delta(\phi))}]=\frac{1}{\sqrt{{\det}_2(I+\phi)}},
\end{equation}
where  ${\det}_2(I+\phi)$ is the Carleman-Fredholm determinant of $I+\phi$.
\end{prop}
The previous proposition remains true if we replace the assumption $||\phi||_{L^2(\R_+^2)}<1$ by the assumption that the eigenvalues of the integral operator with kernel $\phi$ are greater than $-1$ (cf. \cite{chiang} and \cite{grasselli}). 

\begin{prop}\label{prop:15}
Let $F$ be a solution of (\ref{eq:220}) that satisfies the assumptions of Theorem \ref{thm:poppe} and set $\phi_{(x,t)}(a,b):=F(a+b+2x,t)$.  Then for all $x$ and $t$ such that $||\phi_{(x,t)}||_{L^2(\R_+^2)}<1$,
\begin{equation}\label{eq:235}
\tau(x,t):=\E[e^{-\frac{1}{2}\delta(\delta(\phi_{(x,t)}))-\frac{1}{4}\int_0^\infty F(s+2x,t)ds}]^{-2}
\end{equation}
is a tau function for KdV equation (\ref{KdV}) associated to the solution
$$
u(x,t):=4\partial_x^2\log \E[e^{-\frac{1}{2}\delta(\delta(\phi_{(x,t)}))-\frac{1}{4}\int_0^\infty F(s+2x,t)ds}].
$$
\end{prop}

\noindent The assumption $||\phi_{(x,t)}||_{L^2(\R_+^2)}<1$ 
defines a set of $(x,t)$ such that (\ref{eq:235}) is well defined. A similar condition can be found in \cite{tau} and \cite{taniguchi} for the $N$-soliton case. In particular in \cite{taniguchi} the solution is given only for $x\geq 0$.\\

\noindent\textbf{Proof of Proposition \ref{prop:15}.}
The subscripts $(x,t)$ on $\phi$ will be omitted in the following. From Proposition \ref{prop:privault},
$$(\det(I+\phi))^{-1/2}=\E[e^{-\frac{1}{2}\delta(\delta(\phi))-\frac{1}{2}\text{Tr\ }\phi}],$$
where $\text{Tr\ }\phi=\int_0^\infty \phi(s,s)ds=\frac{1}{2}\int_0^\infty F(s+2x,t)ds.$ Then from Theorem \ref{thm:poppe},
$$\tau(x,t):=\det(I+\phi)=\E[e^{-\frac{1}{2}\delta(\delta(\phi_{(x,t)}))-\frac{1}{4}\int_0^\infty F(s+2x,t)ds}]^{-2}$$
is a tau function of KdV and $u=-2\partial_x^2\log \tau$ is solution of KdV.
\hfill \framebox[0.6em]{}\\


\noindent We now consider a particular solution of (\ref{eq:220}) defined by
\begin{equation}\label{eq:250}
F(x,t):=\int_0^\infty e^{8\kappa^3t-\kappa x}d\mu(\kappa),
\end{equation}
where $\mu$ is a $\sigma$-finite measure such that $\mu(0)=0$.
We assume the following holds for some subset $I\subset \R\times \R_+$. 
\begin{equation}\label{eq:236}
\int_0^1 \frac{1}{\kappa}e^{8\kappa^3t-2\kappa x}d\mu(\kappa)<\infty \quad\text{and}\quad\int_1^\infty \frac{1}{\sqrt{\kappa}}e^{8\kappa^3t-2\kappa x}d\mu(\kappa)<\infty,\quad \forall (x,t)\in I.
\end{equation}
We will use the notation $\delta(h)=\int_0^\infty h(s)dW_s$. 

\begin{prop}\label{corol:8}
Let $F$ be given by (\ref{eq:250}) where $\mu$ satisfies (\ref{eq:236}). Let $(X_\kappa)_{\kappa>0}$ be a centered Gaussian process with covariance function $\E [X_{\kappa_1}X_{\kappa_2}]=1/({\kappa_1}+{\kappa_2})$. Then for all $(x,t)\in I$ such that $||\phi_{(x,t)}||_{L^2(\R_+^2)}<1$,
\begin{equation}\label{eq:237}
u(x,t)=4\partial_x^2\log \E\exp\left\{-\frac{1}{2}\int_0^\infty e^{8\kappa^3t-2\kappa x}X_\kappa^2d\mu(\kappa)\right\}
\end{equation}
is solution of KdV equation.
\end{prop}

\noindent We notice the similarity between the covariance function of the process $\left(e^{-\kappa x}X_\kappa\right)_\kappa$ and the matrix (\ref{eq:240}) appearing in the expression of the $N$-soliton solution of KdV (\ref{KdV}).\\

\noindent\textbf{Proof of Proposition \ref{corol:8}.}
Let  $\phi_{(x,t)}(a,b):=F(a+b+2x,t)$. Then
$$\delta(\delta(\phi_{(x,t)}))=\int_0^\infty\int_0^\infty \int_0^\infty e^{8\kappa^3t-\kappa(a+b+2x)}d\mu(\kappa)\ dW_a\ dW_b.$$
By Fubini's theorem for the $\delta$ and $\mu$ applied twice (cf. Theorem \ref{th:fubini} in Appendix),
\begin{align*}
\delta(\delta(\phi_{(x,t)}))&=\int_0^\infty\int_0^\infty \int_0^\infty e^{8\kappa^3t-\kappa(a+b+2x)}dW_a\ dW_b\ d\mu(\kappa)\\
&=\int_0^\infty e^{8\kappa^3t-2\kappa x}\left\{\left(\int_0^\infty e^{-\kappa s}dW_s\right)^2-\int_0^\infty e^{-2\kappa a}ds\right\}d\mu(\kappa)\\
&=\int_0^\infty e^{8\kappa^3t-2\kappa x}X_\kappa^2d\mu(\kappa)-\int_0^\infty \frac{1}{2a}e^{8\kappa^3t-2\kappa x}d\mu(\kappa),
\end{align*}
where $X_\kappa:=\int_0^\infty e^{-\kappa s}dW_s$. The result follows from Proposition \ref{prop:15}.
\hfill \framebox[0.6em]{}\\

The previous proposition gives solutions of KdV different from $N$-solitons and constructed from a Gaussian process with a covariance function corresponding to a particular infinite dimensional extension of a Cauchy matrix. When $\mu$ is a sum of Dirac distributions, we retrieve as a corollary the $N$-soliton solution of KdV.

\begin{corol}\label{corol:7}
Let $c_n$ and $\kappa_n$, $n=1,\ldots,N$ be positive constants. Let $F$ be defined by
\begin{equation}\label{eq:301}
F(x,t):=\sum_{n=1}^Nc_n^2e^{8\kappa_n^3t-\kappa_n x}.
\end{equation}
Then the $N$-soliton solution to KdV equation (\ref{KdV}) is 
\begin{equation}\label{eq:233}
u(x,t)=4\partial_x^2\log \E\exp\left\{-\frac{1}{2}\sum_{n=1}^Nc_n^2e^{8\kappa_n^3t} \left(\int_0^\infty e^{-\kappa_n(x+s)}dW_s\right)^2\right\}.
\end{equation}
\end{corol}

\bigskip
\noindent The solution of KdV obtained in this corollary reminds us of equation (\ref{eq:taniguchi}) proposed by Ikeda where $t$ acts as a parameter and the processes $\left(\int_0^\infty e^{-\kappa_n(x+s)}dW_s\right)_x$ corresponding to the Ornstein-Uhlenbeck processes considered in \cite{taniguchi}. \\

\begin{rmq}
The Corollary \ref{corol:7} can by shown directly from the expectation in the righthand side of (\ref{eq:233}). Indeed, let $\tau(x,t)=\E\exp\left\{-\frac{1}{2}\sum_{n=1}^Nc_n^2e^{8\kappa_n^3t-2\kappa_nx} \left(\int_0^\infty e^{-\kappa_n s}dW_s\right)^2\right\}.$ Let $R:=\text{diag}\{c_n^2e^{8\kappa_n^3t-2\kappa_nx},n=1,\ldots,N \}$ and let $X$ be the $N$-dimensional vector whose components are $X_n:=\int_0^\infty e^{-\kappa_n a}dW_a$. Then $X$ is a Gaussian with mean $0$ and covariance matrix $\Lambda$ given by the Cauchy matrix
$$\Lambda_{m,n}=\E(X_mX_n)=\int_0^\infty e^{-(\kappa_m+\kappa_n)s}ds=\frac{1}{\kappa_m+\kappa_n}.$$
Finally,
$$\tau(x,t)=\E e^{-\frac{1}{2}X^\star R X}=\det(I+R\Lambda)^{-1/2},$$
and $u(x,t)=-2\partial_x^2\log \det(I+R\Lambda).$
Using the identity $\det(I+AB)=\det(I+BA)$ we retrieve the determinant of the matrix (\ref{eq:240}). Hence $u$ is the $N$-soliton solution of KdV (\ref{KdV}).
\end{rmq}

Proposition \ref{prop:15} and Proposition \ref{corol:8} can be extended to the KP hierarchy. Indeed an extension of Theorem \ref{thm:poppe} to KP hierarchy is proved by P{\"o}ppe and Sattinger in \cite{poppe}. Their result is based on a general scheme called dressing method introduced by Zakharov and Shabat \cite{zakharov} for applying the inverse scattering problem method to nonlinear equations. Let us recall the notation for the hierarchy variables ${\bf x}:=(x_1, x_2,\ldots)$ and ${\bf z}:=(z_1, x_2,\ldots)$. Then P{\"o}ppe and Sattinger showed that a tau function for KP hierarchy can be written as the Fredholm determinant of the integral operator
\begin{equation}\label{eq:241}
\mathcal{F}_{\bf{x}}\psi(y)=\int_0^\infty F((x_1+y,x_2,x_3,\ldots),(x_1+z,x_2,x_3,\ldots))\psi(z)dz,
\end{equation}
where $F({\bf x},{\bf z})$ satisfies the following system of linear pde in infinitely many variables
\begin{equation}\label{eq:242}
\frac{\partial}{\partial x_n}F-\frac{\partial^n}{\partial x_1^n}F+(-1)^n\frac{\partial^n}{\partial z_1^n}F=0,\quad n=2,3,\ldots,
\end{equation}
and $\mathcal{F}_{\bf{x}}$ is assumed to belong to $LC_\nu$, $1/2 <\nu\leq 1$ for all ${\bf{x}}$.\\

\noindent System (\ref{eq:242}) is the analogous of (\ref{eq:220}) for KP hierarchy. The variables $x_2, x_3,\ldots$ in equation (\ref{eq:241}) act as parameters exactly like the time $t$ does in equation (\ref{eq:235}). Defining $\xi({\bf x},k):=\sum_{j=1}^\infty x_jk^j$, system (\ref{eq:242}) admits the solution
$$F({\bf x},{\bf z})=e^{\xi({\bf x},p)-\xi({\bf z},q)}.$$
More general solutions can then be obtained by superposition of such fundamental solutions as follows
\begin{equation}
\label{eq:150}
F({\bf x},{\bf z})=\int_{\C^2}e^{\xi({\bf x},p)-\xi({\bf z},q)}d\mu(p,q),
\end{equation}
where $\mu$ is some measure on $\C^2$. Its support is included in the set $\{(p,q)\in\C^2, \text{Re\ }p<0<\text{Re\ }q\}$.\\ 

\noindent In particular, solutions to KP equation (\ref{eq:KP}) are expressed from the tau function $\tau$ of KP hierarchy as
$$u(x,y,t)=-2\partial_x^2 \log \tau(x,y,t,0,\ldots),$$
for all $(x,y,t)$ such that $\tau(x,y,t,0,\ldots)$ is positive (cf. \cite{poppe}).

\begin{prop}\label{prop:16}
Let $1/2 <\nu\leq 1$ and let $\mathcal{F}_{\bf{x}}\in LC_\nu$ be an integral operator defined by (\ref{eq:241}) where $F$ is real valued and satisfies (\ref{eq:242}). Let $\phi_{\bf x}(y,z):=F((x_1+y,x_2,x_3,\ldots),(x_1+z,x_2,x_3,\ldots))$. Then for all ${\bf x}$ such that $||\phi_{\bf x}||_{L^2(\R_+^2)}<1$,
\begin{equation}\label{eq:245}
\tau({\bf x}):=\E[e^{-\frac{1}{2}\delta(\delta(\phi_{\bf x}))-\frac{1}{2}\int_0^\infty \phi_{\bf x}(y,y)dy}]^{-2}
\end{equation}
is a tau function for KP hierarchy.
\end{prop}

\noindent\textbf{Proof of Proposition \ref{prop:16}.}
The proof of this proposition follows the same lines as the proof of Proposition \ref{prop:15}. We leave the details to the reader.
\hfill \framebox[0.6em]{}\\

\section{Appendix}

\subsection{Fubini's theorem for Skorohod integrals}

A Fubini's theorem for Skorohod integrals can be found in \cite{nualart} (Ex. 3.2.7). It is not difficult to prove the following version where the integration domain is not necessarily bounded.

\begin{thm}\label{th:fubini} Consider a random field $\{u_t(x),0\leq t<\infty, x\in G\}$, where $G\subset\R$, such that for all $x\in G$, $u_\cdot(x)\in \text{Dom\ }\delta$. Let $\mu$ be a $\sigma$-finite measure on $G$. Suppose that $\int_G \left(\E\int_0^\infty|u_t(x)|^2dt\right)^{1/2}d\mu(x)<\infty$, and $\int_G \E[\delta(u_\cdot(x))^2]^{1/2}d\mu(x)<\infty,$ then the process $\{\int_Gu_t(x)d\mu(x),0\leq t<\infty\}$ is Skorohod integrable and
$$\delta\left(\int_Gu_t(x)d\mu(x)\right)=\int_G\delta(u_\cdot(x))d\mu(x).$$
\end{thm}

\newpage
\renewcommand{\refname}{Bibliography}
\bibliographystyle{plain}
\bibliography{bibliographie}

\end{document}